\documentclass{micls}%ajcnosec}
\usepackage{color,graphicx}
\usepackage{amssymb}
\usepackage{hyperref}

\newauthor{%
Italo J. Dejter, Luis R. Fuentes}{%
Italo J. Dejter, Luis R. Fuentes}{%
University of Puerto Rico\\
Rio Piedras, PR 00936-8377}[%
italo.dejter@gmail.com, luis.fuentes@upr.edu]

\newauthor{%
Carlos A. Araujo}{%
Carlos A. Araujo}{%
Universidad del Atl\'antico\\
Barranquilla, Colombia}[%
carlosaraujo@mail.uniatlantico.edu.co]

\title{There is but one PDS in $\mathbb{Z}^3$ inducing just square components}%{Uniqueness of Perfect Dominating Set in $\mathbb{Z}^{3}$ Inducing 4-Cycles}

\keywords{perfect dominating sets; unit distance graph; integer lattice}

\classnbr{Primary 05C69; Secondary 94B25}

\begin{document}
\begin{abstract} It is known that in the unit distance graph of the lattice $\mathbb{Z}^3\subset\mathbb{R}^3$ there exists a dominating set $S$ with $4$-cycles as sole induced components and each vertex of $\mathbb{Z}^3\setminus S$ having a unique neighbor in $S$. We show $S$ is unique.
\end{abstract}

\section{PERFECT DOMINATING SETS, (PDS\thinspace s)}\label{s1}

Let $\Gamma=(V,E)$ be a graph and let $S\subset V$.
The closed neighborhood of a vertex $\theta\in V$ in $\Gamma$ is denoted $N[\theta]$.
Let $[S]$ be the subgraph of $\Gamma$ induced by $S$.
The induced components of $S$, namely the connected components of $[S]$ in $\Gamma$, are said to be the {\it components} of $S$.
Several definitions of perfect dominating sets in graphs are considered in the literature \cite{HHS,Klos}. We work with the following one \cite{W} denoted with the short acronym PDS, to make a distinctive difference:

\vspace*{1mm}

$S$ is a PDS of $\Gamma$ $\Leftrightarrow$ each vertex of $V\setminus S$ has a unique neighbor in $S$.

\vspace*{1mm}

\noindent This definition (of PDS) differs from that of a
`perfect dominating set' as in \cite{FH,GS,LS} (that for us is a stable PDS coinciding with the perfect code of \cite{Biggs} or with the efficient dominating set of \cite{BBS,DS,HHS}). With our not necessarily stable definition of perfect dominating set, denoted PDS, our main result, stated below as Theorem~\ref{mr}, has a narrowing spirit as that of Theorem 2.6 of just cited \cite{LS}.

Let $0<n\in\mathbb{Z}$. The following graphs are considered. The unit distance graph $\Lambda_n$ of the $n\cite{}$-dimensional integer lattice $\mathbb{Z}^n\subset\mathbb{R}^n$ has vertex set $\mathbb{Z}^n$ and exactly one edge between
each two vertices if and only if their Euclidean distance is 1.
An $n$-{\it cube} is the cartesian graph product $Q_n=K_2\square K_2\square\cdots\square K_2$ of precisely $n$ copies of the complete graph $K_2$. In particular, a 2-cube $Q_2$ is a square, that is a 4-cycle.
A {\it grid graph} is the cartesian graph product of two path graphs.

Our definition of a PDS $S$ allows components of $S$ in $\Gamma$ which are not isolated vertices. For example:
{\bf(a)} tilings with generalized Lee $r$-spheres, for fixed $r$ with $1<r\le n$ in $\mathbb{Z}$ (e.g., crosses with arms of length one if $r=n$), furnish $\Lambda_n$ with PDS\thinspace s whose components are $r$-cubes \cite{Etzion}, including that of our Theorem~\ref{mr}, below; (It is most remarkable that $r=n$ $\Leftrightarrow$ $n\in\{2^r-1,3^r-1;0<r\in\mathbb{Z}\}$ \cite{BE}); {\bf(b)} {\it total perfect codes} \cite{AD,KG}, that is PDS\thinspace s whose components are copies of $K_2=P_2$ in the $\Lambda_n$\thinspace s and grid graphs; (these appear as {\it diameter perfect Lee codes} \cite{E,BH});
{\bf(c)} PDS\thinspace s in $n$-cubes \cite{BD,D,DD,DP,DW,W}, where $0<n\in\mathbb{Z}$, including the perfect codes of \cite{DM}; {\bf(d)} PDS\thinspace s in grid graphs \cite{DD,KG}.

\begin{theorem}\label{mr} There is only one PDS in $\Lambda_3$ whose components are 4-cycles.
\end{theorem}

This is proved as Theorem~\ref{t3} once some auxiliary notions are presented.

\section{INDUCED COMPONENTS}\label{s2}

The {\it distance} $d(u,v)$ between two vertices $u$ and $v$ of $\Lambda_n$
is defined as the minimum length of any path connecting $u$ and
$v$. The following is an elementary extension of a result of
\cite{W} for $n$-cubes.

\begin{theorem}\label{t1}
Let $S$ be a {\rm PDS} in $\Lambda_n$. Let $J_S$ be a set of indices $j$ for the corresponding components $S^j$ of $S$.
Each $S^j$ is a cartesian graph product of connected subgraphs of
$\Lambda_1$. Thus, if such $S^j$ is a finite subgraph $\Theta$ of $\Lambda_n$, then $S^j$ is of the form $P_{i_1^j}\square P_{i_2^j}\square\cdots\square P_{i_n^j},$ where $P_{i_k^j}$ is a path of length $i_k^j-1\ge 0$, for $k=1,\ldots,n$.
\end{theorem}

A {\rm PDS} in $\Lambda_n$ whose components are all
isomorphic to a fixed finite graph $\Theta$ (as in Theorem~\ref{t1}) is called a {\rm
PDS}$[\Theta]$.
If no confusion arises, $n$-tuples representing elements
of $\mathbb{Z}^{n}$ are written with neither commas nor external parentheses. We denote $00\ldots 0=O$, $10\ldots 0=e_1$, $010\ldots 0=e_2$,
$\ldots $, $00\ldots 1=e_n$.

At the end of Section 6 of \cite{Etzion} (in the original setting of item (a) above in Section~\ref{s1}), all the indices $i_k^j$ of our Theorem~\ref{t1} are shown to be less than 2.

\section{LATTICE-LIKE DOMINATING SETS}\label{s3}

Let $\Theta=(V,E)$ be a finite subgraph of $\Lambda_n$ and let $z\in\mathbb{Z}^n$.
Then $\Theta+z$ denotes the graph $\Theta'=(V',E')$, where
$V'=V+z=\{w$; there is $v\in V,w=v+z\}$ and $uv\in E\Leftrightarrow(u+z)(v+z)\in E'$. Let $S$ be a PDS$[\Theta]$ and let a copy $D$ of $\Theta$ be a component of $S$. Then $S$ is
said to be \textit{lattice-like} if there is a lattice $L$ (that is, a subgroup $L$ of $\mathbb{Z}^n$) so that $D'$ is a component of $S$ if and only if there is $z\in L$ with $D'=D+z$. Examples above (\cite{Etzion,BE,E,BH}) are lattice-like.

If $S$ is a PDS$[\Theta]$ with $\Theta=(V,E)$, then $S$ can be seen as a tiling of $\mathbb{Z}^n$ by the induced subgraph $\Theta^*$ of $\Lambda_n$ on the set $V^*=\{v\in\mathbb{Z}^n;d(v,V)\leq 1\}$. Thus, a lattice-like tiling will be understood in the same way as a lattice-like PDS.
We need the following form of Theorem 6 \cite{BH} for the proof of Theorem~\ref{t3}. Recall that given a graph $G$, the distance $d(v,H)$ between a vertex $v$ of $G$ and a subgraph $H$ of $G$ is the shortest distance between $v$ and the vertices of $H$.

\begin{theorem}\label{C} Let $\Theta$ be a subgraph of $\Lambda_n$. Let $\Theta^*$ be an induced supergraph of $\Theta$ in $\Lambda_n$ such that a vertex $v$
is in $\Theta^*$ if and only if $d(v,\Theta)\leq 1$. Let $D=(V,E)$ be a
copy of $\Theta^*$ that contains vertices $O,e_1,\ldots ,e_n$. Then, there is a lattice-like
PDS$[\Theta]$ if and only if there exists an abelian group $G$ of order $|V|$
and a group epimorphism $\Phi :\mathbb{Z}^n\rightarrow G$ such that the
restriction of $\Phi$ to $V$ is a bijection.
\end{theorem}

\section{THE PROOF}

\begin{theorem}\label{t3}
There do not exist non-lattice-like PDS$[Q_2]$\thinspace s in $\Lambda_3$. In addition,
there exists exactly one lattice-like PDS$[Q_2]$ in $\Lambda_3$.
\end{theorem}

\begin{figure}[htp]\hspace*{1mm}%Figure 1
  \includegraphics[scale=0.35]{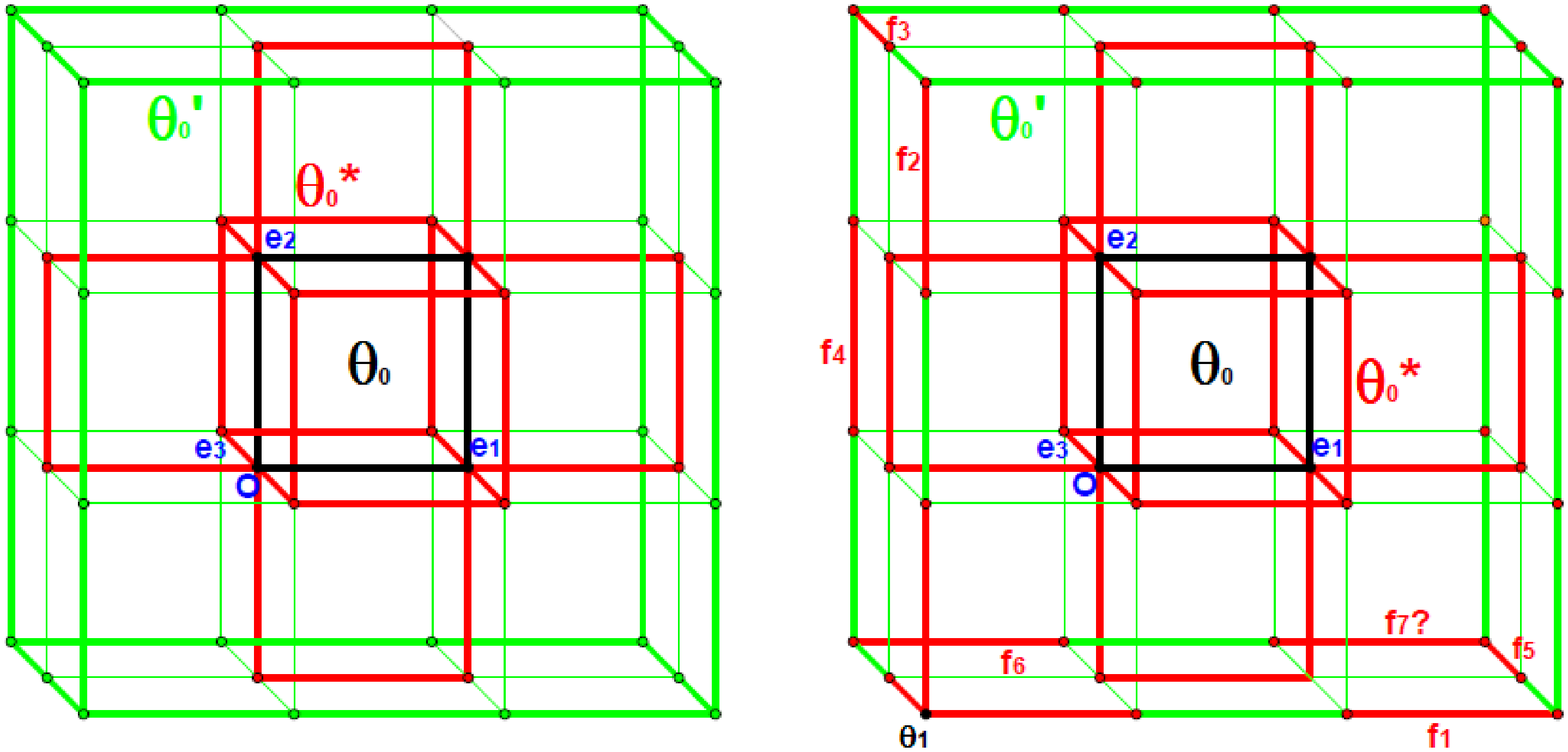}
\caption{$\Theta_0\subset\Theta_0^*\subset\Theta_0'$, and the case of one corner of $\Theta_0'-\Theta_0^*$ in $S$}%\label{}
\end{figure}

\begin{figure}[htp]\hspace*{18mm}%Figure 2
  \includegraphics[scale=0.50]{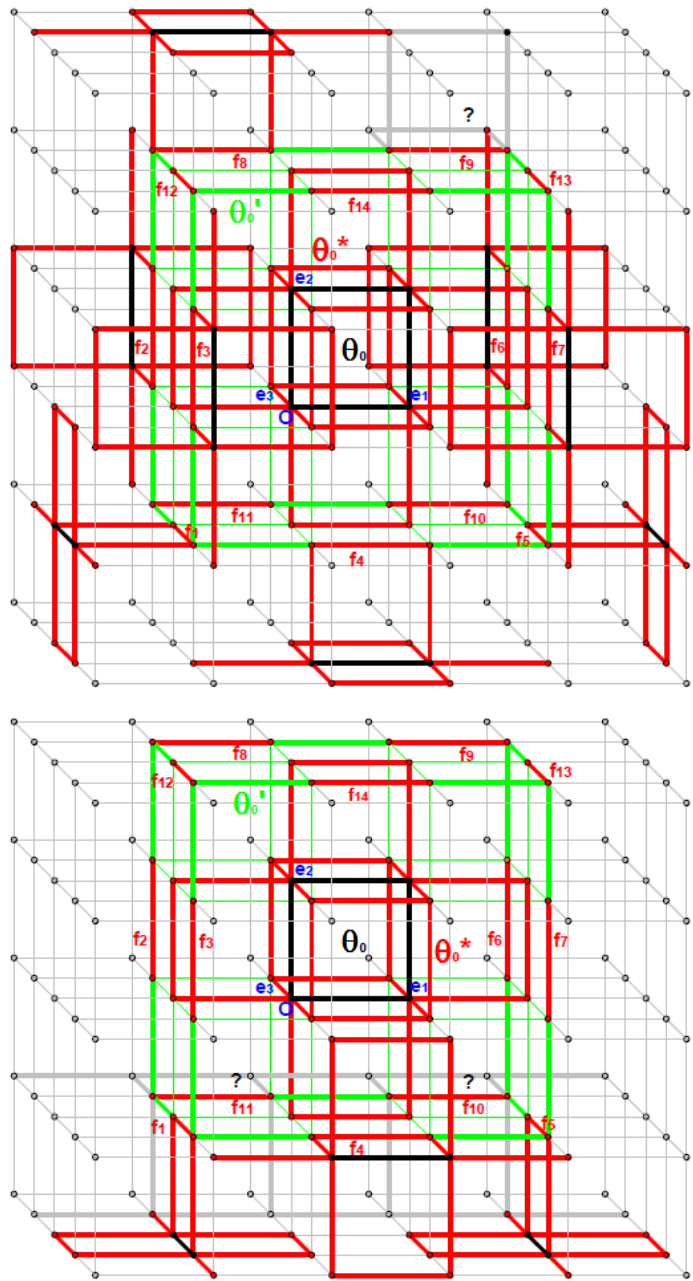}
\caption{Instances of no corners of $\Theta_0'-\Theta_0^*$ in $S$, (a$_1$) and (a$_2$)}
\end{figure}

\begin{proof} Theorem 8 \cite{Etzion} insures the
existence of a PDS$[Q_2]$ in $\Lambda_3$. In fact, the connected
components of such PDS$[Q_2]$ are the generalized Lee spheres
$S_{3,2,0}$ inside the corresponding generalized Lee spheres
$S_{3,2,1}$ (in their inductive construction in Section 1
\cite{Etzion}) that form the lattice tiling $\Lambda_{3,2}$ (in the
notation of \cite{Etzion}) insured by that Theorem 8. According to
the theorem, this $\Lambda_{3,2}$ has generator matrix (as defined
in Section 3 \cite{Etzion}):

\begin{eqnarray}\label{gm}
\left(
          \begin{array}{ccc}
            1 & 0 & 3 \\
            0 & 2 & 5 \\
            0 & 0 & 10 \\
          \end{array}
        \right)
\end{eqnarray}

In terms of Theorem~\ref{C}, the generator matrix~(\ref{gm}) corresponds to the group epimorphism
$\Phi:\mathbb{Z}^3\rightarrow G=\mathbb{Z}_{20}$ given by $\Phi(e_1)=2$; $\Phi(e_2)=5$ and $\Phi(e_3)=6$ where $\Phi$ is obtained first by multiplying the matrix~(\ref{gm}) by an unknown vector and then solving the corresponding system of equations mod 20.
To see that this is the only PDS$[Q_2]$ in $\Lambda_3$, we note that there are only two possible abelian groups $G$ for the epimorphism $\Phi$, namely: $G=\mathbb{Z}_{20}$ and $G=\mathbb{Z}_2\times \mathbb{Z}_2 \times \mathbb{Z}_5$. It can be easily checked \cite{Lucho} that {\bf(a)} there are
just $32$ epimorphism from $\mathbb{Z}^3$ onto $G=\mathbb{Z}_{20}$ and none from $\mathbb{Z}^3$ onto $G=\mathbb{Z}_{2}\times \mathbb{Z}_2\times \mathbb{Z}_5$; {\bf(b)} every possible assignment for $\Phi(e_1)$, $\Phi(e_2)$ and $\Phi(e_3)$ has order $10$, $10$ and $4$, respectively, in $\mathbb{Z}_{20}$. As a result, all 4-cycles induced by each lattice-like PDS$[Q_2]$ associated (via Theorem~\ref{C}) to a corresponding of these $32$ epimorphisms are placed in the same way in $\Lambda_3$. Each such lattice-like PDS$[Q_2]$ in $\Lambda_3$ is equivalent to the one obtained via matrix~(\ref{gm}).

Assume there is a non-lattice-like PDS$[Q_2]$ $S$ in $\Lambda_3$ so that the components of $[S]$ are 4-cycles $Q_2$; let $\Theta=\Theta_0$ be such a component. We may assume that $\Theta_0$ has vertices $O$, $e_1$, $e_2$, $e_1+e_2$. The graph $\Theta^*=\Theta_0^*$ is contained in a graph $\Theta_0'$ isomorphic to $P_4\square P_4\square P_3$ as on the left of Figure 1, where $\Theta_0$ has its edges thick black, the rest of $\Theta_0^*$ has them  red and the rest of $\Theta_0'$ has them green, thick for the paths between the eight corners (vertices of degree 3 in $\Theta_0'$: $-e_1-e_2\pm e_3$, $2e_1-e_2\pm e_3$, $-e_1+2e_2\pm e_3$, $2e_1+2e_2\pm e_3$) and thin for the rest. The realization of $\Theta_0^*$ in $\mathbb{R}^3$ has convex hull containing tightly  $\Theta_0'$. Similar colors and traces are used in the representations in Figures 2-14, where:
{\bf(I)} The red thin-trace lines un Figures 2-3 and 4-5 represent edges incident to vertices in subgraphs $\Theta^*$ (red thick-trace edges) involved in our arguments by contradiction, indicated by question marks (?);
{\bf(II)} yellow squares in Figures 7-14 indicate where to paste accordingly (the front of) the top and (the back of) the lower parts of each figure to obtain a continuation of the lattice representation in each figure; {\bf(III)} edges not mentioned in (I) or (II) are traced in dashed green color.

Assume no vertex of $\Theta_0'-\Theta_0^*$ is in $S$.
By symmetry there is a 1-factor $F$ in $\Theta_0'-\Theta_0^*$ each of whose edges has an endvertex $\notin V(\Theta_0'-\Theta_0^*)$ dominated by a vertex in a 4-cycle induced by $S$. In each case we will reach a contradiction: $F$ is either as in case (a) or (b) below, depending on the feasible dispositions of four edges of $F$ over the four maximal paths of length 2 between the eight corners of $\Theta_0'$, namely either with their eight endvertices having convex hull tightly containing a copy of $P_4\square P_4\square P_2$ (say convex hull $[-1,2]\times[-1,2]\times[-1,0]$) or not (in which case partial convex hulls $\{-1\}\times[-1,2]\times[-1,0]$ and $\{ 2\}\times[-1,2]\times[ 0,1]$ appear, not leading to a total convex hull as above), that we have respectively either as the four edges $f_1,f_5,f_{12},f_{13}$, for (a), or as the four edges $f_1,f_4,f_8,f_{12}$ for (b). These instances are: (with (a) further subdivided into subcases (a$_1$) and (a$_2$), below)

{\bf(a)} (Figure 2, top) The edges of $F$ are:
$$\begin{array}{lll}
^{f_1\,\,=(-e_1-e_2,-e_1-e_2-e_3),}_{f_4\,\,=(-e_2-e_3,\,\,e_1-e_2-e_3),}&
^{f_2\,\,=(-e_1+e_3,-e_1+e_2+e_3),}_{f_5\,\,=(\,2e_1-e_2,\,\,2e_1-e_2-e_3),}&
^{f_3\,\,=(-e_1-e_3,-e_1+e_2-e_3),}_{f_6\,\,=(2e_1+e_3\,,2e_1+e_2+e_3),}\\
^{f_7\,\,=(\,2e_1-e_3,\,\,2e_1+e_2-e_3),}_{f_{10}=(e_1-e_2+e_3,2e_1-e_2+e_3),}&
^{f_8\,\,=(2e_2+e_3,-e_1+2e_2+e_3),}_{f_{11}=(-e_2+e_3,\,\,-e_1-e_2+e_3),}&
^{f_9\,\,=(e_1+2e_2+e_3,2e_1+2e_2+e_3),}_{f_{12}=(-e_1+2e_2,-e_1+2e_2-e_3),}\\
^{f_{13}=(\,2e_1+2e_2,2e_1+2e_2-e_3),}&^{f_{14}=(2e_2-e_3,\,\,\,\,\,\,e_1+2e_2-e_3).}&
\end{array}$$

\noindent We may take step by step either option (a$_1$) or option (a$_2$) below (where, instead of saying that a vertex $v$ is dominated by an endvertex of an edge $f$, we simply say that $v$ is dominated by $f$, or that $v\in(f)$), with $(f)$ representing the set of vertices dominated by the endvertices of $f$):

{\bf(a$_1$)} The first eight edges in (a) have each an endvertex
dominated by a vertex in a 4-cycle. The involved 4-cycles contain the following edges:

$\bullet$ $f_1-e_1$ (the translation of $f_1$ via the vector $-e_1$),

$\bullet$ $f_2+e_3$ (forced, since $f_2-e_1$ dominates $-2e_1\in(f_1-e_1)$),

$\bullet$ $f_3-e_3$ (forced, since $f_3-e_1$ contains $-2e_1-e_3\in(f_1-e_1)$),

$\bullet$ $f_4-e_2$ (forced, since $f_4-e_3$ dominates $-e_1-e_2-2e_3\in(f_3-e_3)$),

$\bullet$ $f_5+e_1$ (forced, since $f_4-e_2$ contains $2e_1-2e_2-e_3\in(f_4-e_2)$),

$\bullet$ $f_6+e_3$ (forced, since $f_6+e_1$ dominates  $3e_1\in(f_5+e_1)$),

$\bullet$ $f_7-e_3$ (forced, since $f_7+e_1$ contains  $3e_1-e_3\in(f_5+e_1)$) and

$\bullet$ $f_8+e_2$ (forced, since $f_8+e_3$ dominates $-e_1+2e_2+2e_3\in(f_2+e_3)$).

\noindent Now, there is no way for the edge $f_9$ to be dominated by a copy of $K_2$ external to $\Theta_0'-\Theta_0^*$ (since $f_9+e_2$ contains $e_1+3e_2+e_3\in(f_8+e_2)$ while $f_9+e_3$ contains $2e_1+2e_2+2e_3\in(f_6+e_3)$), a contradiction.

{\bf(a$_2$)} (Figure 2, bottom) The edges $f_1$, $f_5$ and $f_4$ have each one endvertex dominated by a vertex in a 4-cycle containing the respective edges
$f_1-e_2$, $f_5-e_2$ (edge pair not contemplated in case (a$_1$)) and $f_4-e_3$ (forced, since $f_4-e_2$ contains vertex $-2e_2\in(f_1-e_2)$). But then
only one of $f_{11}$ and $f_{10}$ must be dominated by $f_{11}-e_2$ or $f_{10}-e_2$, while the remaining one must be dominated by $f_{11}-e_3$ or $f_{10}-e_3$, which produces a contradiction since $f_{11}-e_2\in(f_1-e_2)$, and $f_{10}-e_2\in(f_5-e_2)$.

\begin{figure}[htp]%\hspace*{37mm}%Figure 3
  \includegraphics[scale=0.30]{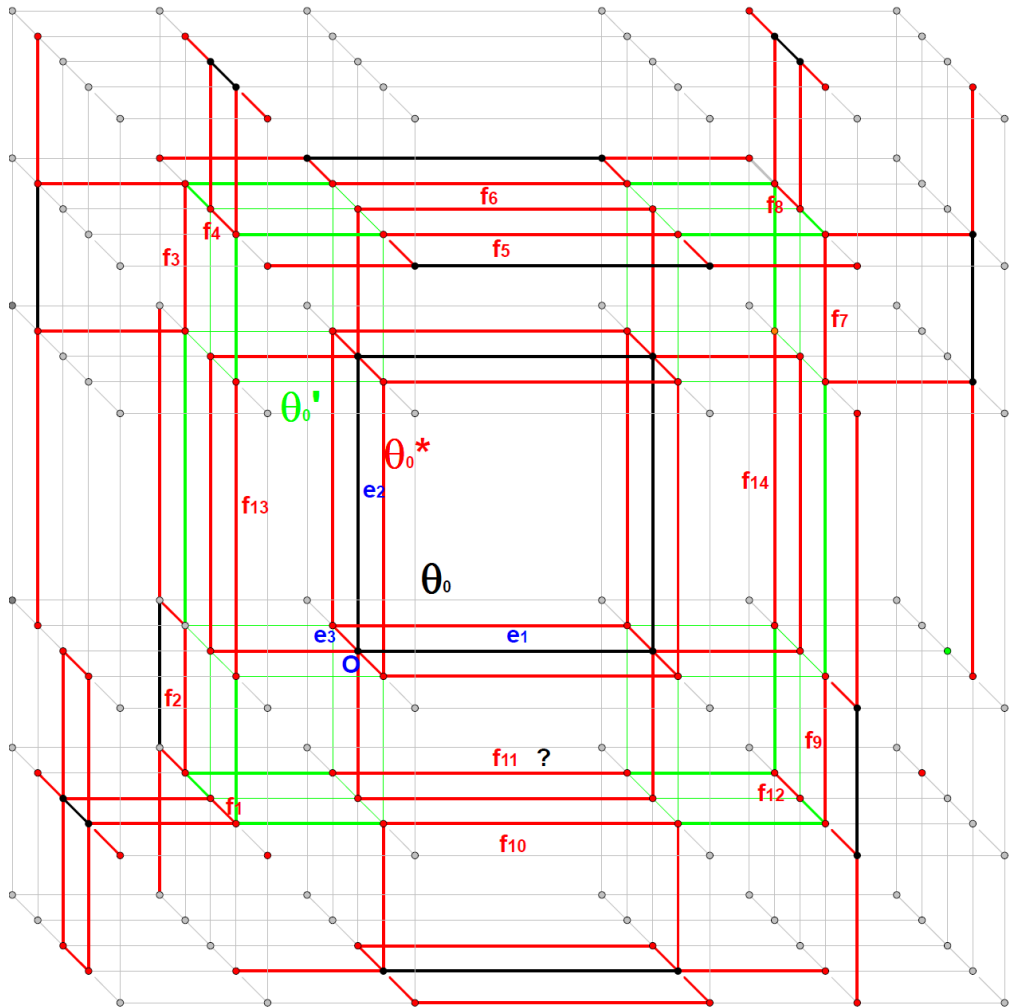}
\caption{Instance of no corners of $\Theta_0'-\Theta_0^*$ in $S$, case (b)}
\end{figure}

\begin{figure}[htp]\hspace*{-5mm}%Figure 4
  \includegraphics[scale=0.40]{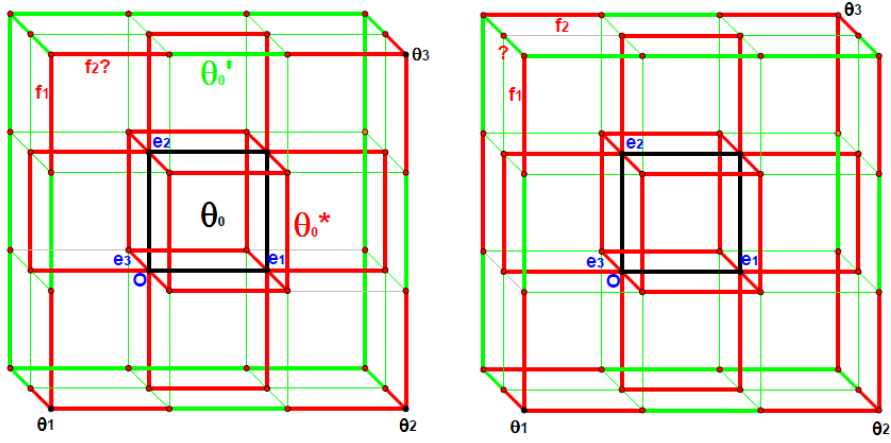}
\caption{The two cases of three corners of $\Theta_0'-\Theta_0^*$ in $S$}
\end{figure}

{\bf(b)} (Figure 3)
The edges of $F$ are:
$$\begin{array}{lll}
^{f_1\,\,=(-e_1-e_2\,\,,\,-e_1-e_2-e_3),}_{f_4\,\,=(-e_1+2e_2,-e_1+2e_2-e_3),}&
^{f_2\,\,=(-e_1+e_3\,\,,\,-e_1-e_2+e_3),}_{f_5\,\,=(\,2e_2-e_3\,\,,\,\,\,e_1+2e_2-e_3),}&
^{f_3\,\,=(-e_1+e_2+e_3,-e_1+2e_2+e_3),}_{f_6\,\,=(\,\,2e_2+e_3\,\,,\,\,\,e_1+2e_2+e_3),}\\
^{f_7\,\,=(\,2e_1+2e_2-e_3,\,2e_1+e_2-e_3),}_{f_{10}=(-e_2-e_3\,\,,\,\,\,e_1-e_2-e_3),}&
^{f_8\,\,=(\,\,2e_1+2e_2,2e_1+2e_2+e_3),}_{f_{11}=(-e_2+e_3\,\,,\,e_1-e_2+e_3),}&
^{f_9\,\,=(\,2e_1-e_3\,\,,\,2e_1-e_2-e_3),}_{f_{12}=(\,\,2e_1-e_2\,\,,\,2e_1-e_2+e_3),}\\
^{f_{13}=(-e_1-e_3\,\,,-e_1+e_2-e_3),}&
^{f_{14}=(\,2e_1+e_3\,\,,\,2e_1+e_2+e_3).}&\\
\end{array}$$

We may assume step by step that the first ten edges of $F$ have each an endvertex dominated by the copy of $K_2$ containing respectively:

$\bullet$ $f_1-e_1$,

$\bullet$ $f_2+e_3$ (forced, since $f_2-e_1$ contains $e_3-2e_1\in(f_1-e_1)$),

$\bullet$ $f_3-e_1$ (forced, since $f_3+e_3$ contains $e_2+2e_3-e_1\in(f_2+e_3)$),

$\bullet$ $f_4+e_2$ (forced, since $f_4-e_1$ contains $2e_2-2e_1\in(f_3+e_1)$),

$\bullet$ $f_5-e_3$ (forced, since $f_5+e_2$ contains $3e_2-e_3\in(f_4+e_2)$),

$\bullet$ $f_6+e_3$ (forced, since $f_6+e_2$ contains $3e_2+e_3\in(f_4+e_2)$),

$\bullet$ $f_7+e_1$ (forced, since $f_7-e_3$ contains $2e_1+2e_2-2_3\in(f_5-e_3)$),

$\bullet$ $f_8+e_2$ (forced, since $f_8+e_1$ contains $3e_1+2e_2\in(f_7+e_1)$),

$\bullet$ $f_9-e_3$ (forced, since $f_9+e_1$ contains $3e_1-e_3\in(f_7+e_1)$) and

$\bullet$ $f_{10}-e_2$ (forced, since $f_{10}-e_3$ contains $e_1-e_2-2e_3\in(f_9-e_3)$).

\noindent Now, $f_{11}$ does not have an endvertex dominated by any copy of $K_2$ in the presence of the previous forced dominations of copies of $K_2$ (since $f_{11}-e_2$ dominates $\{-2e_2,-2e_2\}\subset(f_{10}-e_2)$ while $f_{11}+e_3$ contains $2e_3-e_2\in(f_2+e_3)$).

If just one or three corners of $\Theta_0'$ (in this second case, for corner distance triple either $(3,3,6)$ or $(3,5,8)$) were in $S$, the remaining vertices of $\Theta_0'-\Theta_0^*$ forms no 1-factor $F$, contradicting the existence of $S$. (Figure 1, right, and Figure 4).
In the case of one corner, let this corner be $\theta_1=-e_1-e_2-e_3$, which dominates $\theta_1+e_1,\theta_1+e_2$ and $\theta_1+e_3$. Then $F$ must contain:
$$\begin{array}{lll}
^{f_1=(e_1-e_2-e_3,2e_1-e_2-e_3),}_{f_4=(-e_1+e_3,-e_1+e_2+e_3),}&
^{f_2=(-e_1+e_2-e_3,-e_2+2e_2-e_3),}_{f_5=(2e_1-e_2,2e_1-e_2+e_3),}&
^{f_3=(-e_2+2e_2,-e_2+2e_2+e_3),}_{f_6=(-e_1-e_2+e_3,-e_2+e_3).}
\end{array}$$

\noindent Now, $F$ should also contain $f_7=(e_1-e_2+e_3,2e_1-e_2+e_3)$, with its terminal vertex already present in $f_5$, a contradiction. With three corners and distance triple $(3,3,6)$, let these corners be $\theta_1$, $\theta_2=2e_1-e_2-e_3$ and $\theta_3=2e_1+2e_2-e_3$. Then $F$ must contain $f_1=(-e_1+e_2-e_3,-e_1+2e_2-e_3)$ and $f_2=(-e_1+2e_2-e_3,2e_2-e_3)$ that have a vertex in common, a contradiction. With distance triple $(3,5,8)$, let the three corners be $\theta_1$, $\theta_2$ and $\theta'_3=\theta_3+2e_3$. Then $F$ must contain $f_1$ as above and $f'_2=f_2+2e_3$, leaving vertex $-e_1+2e_2$ not in $F$, another contradiction.

\begin{figure}[htp]\hspace*{22mm}%Figure 5
  \includegraphics[scale=0.45]{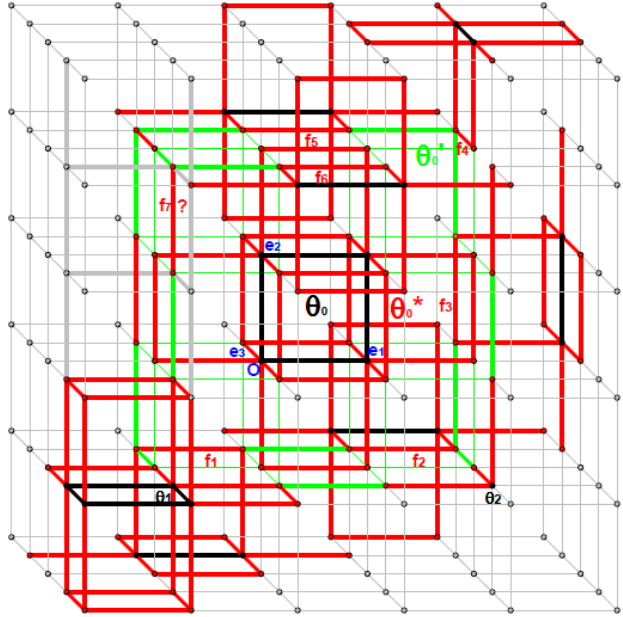}
\caption{The case of two corners of $\Theta_0'-\Theta_0^*$ in $S$ at distance 3}
\end{figure}

\begin{figure}[htp]\hspace*{22mm}%Figure 6
  \includegraphics[scale=0.45]{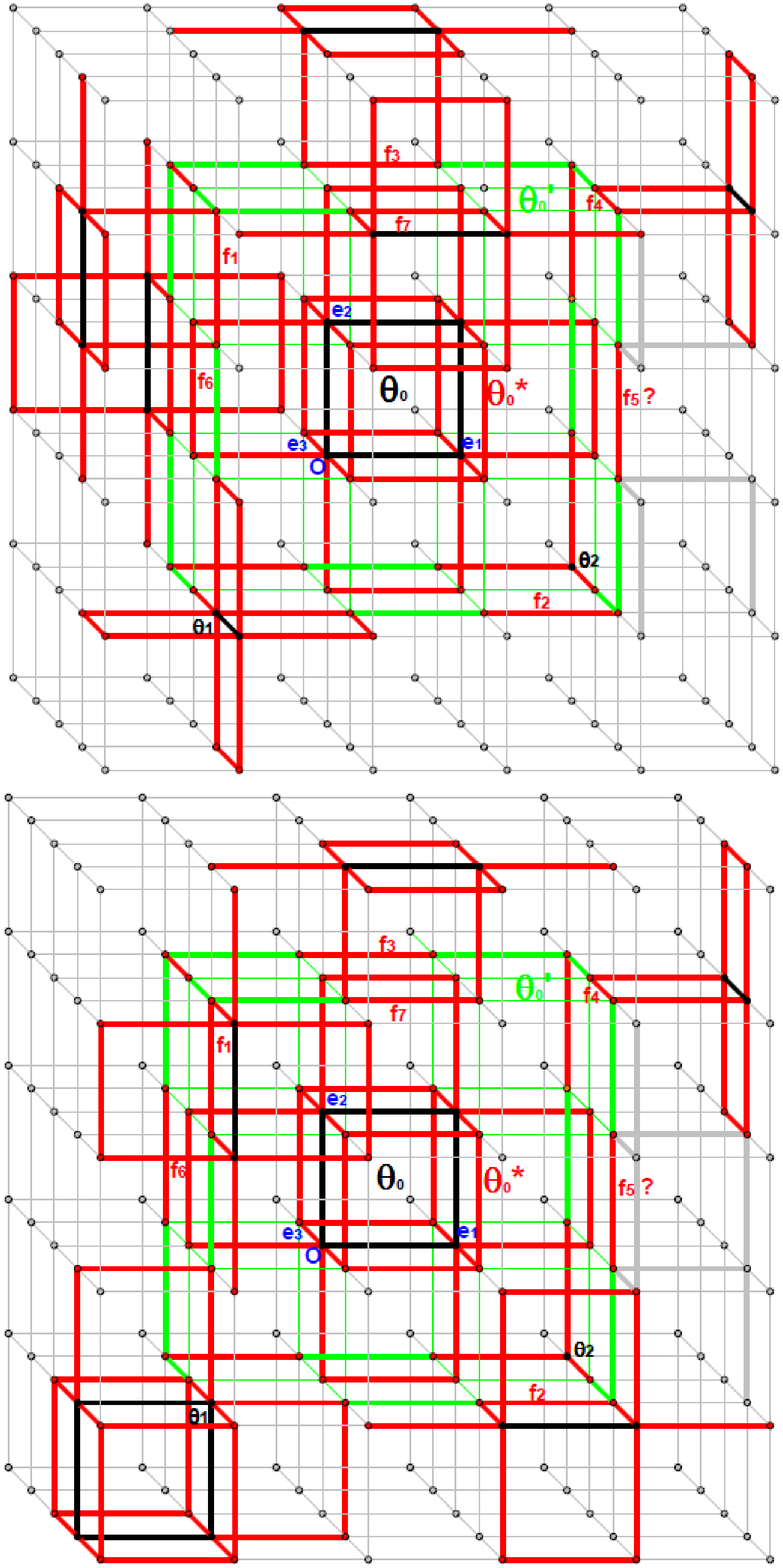}
\caption{The two cases of two corners of $\Theta_0'-\Theta_0^*$ in $S$ at distance 5}
\end{figure}

\begin{figure}[htp]\hspace*{20mm}%Figure 7
  \includegraphics[scale=0.45]{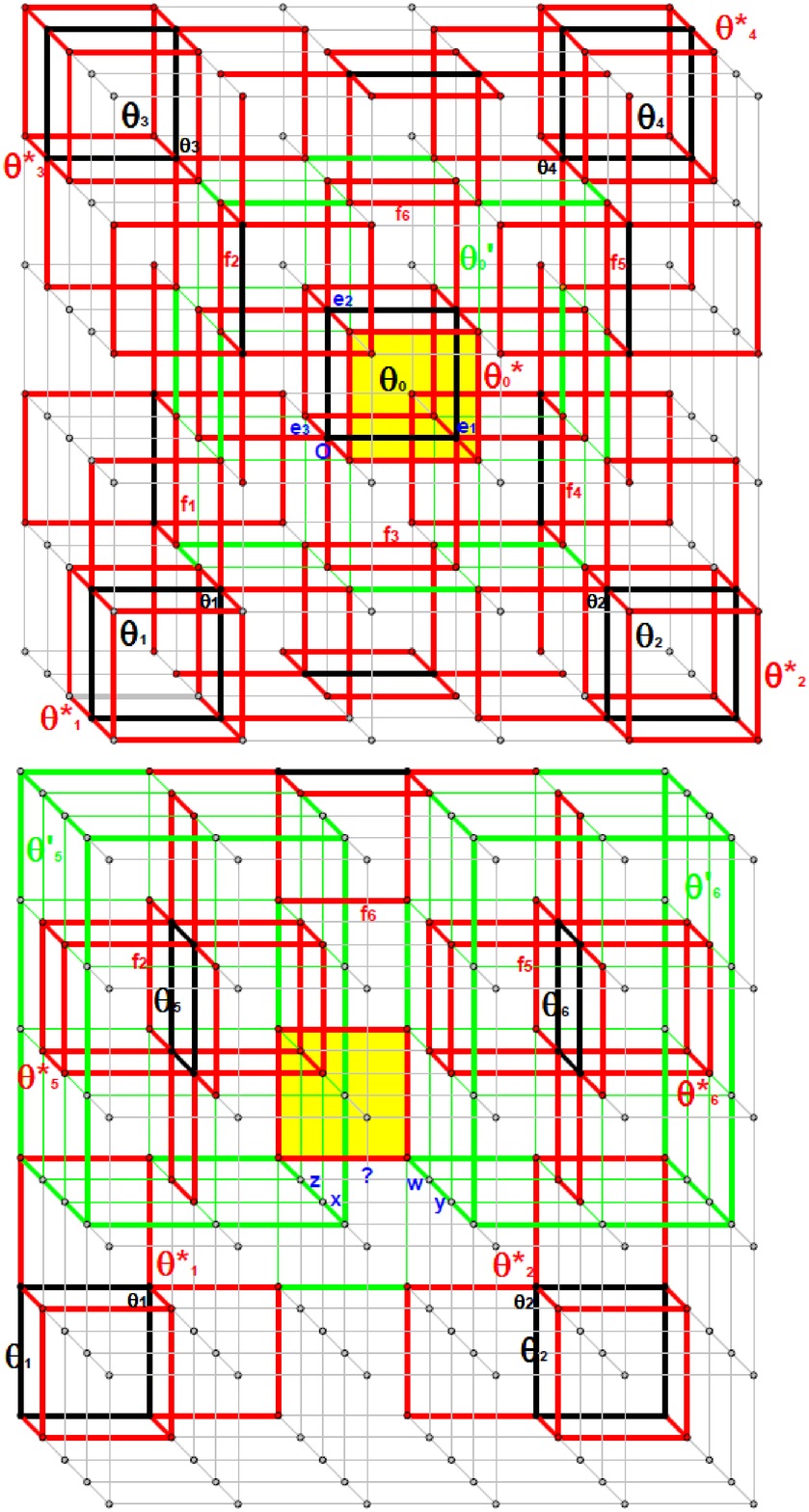}
\caption{Four corners, instance (A), case (a)}
\end{figure}

We will rule out the cases of only two corners of $\Theta_0'$ being in $S$.
If the two are at distance 3 (Figure 5) they may be taken up to symmetry as $\theta_1=-e_1-e_2-e_3$ and $\theta_2=2e_1-e_2-e_3$. In $\Theta_0'-\Theta_0^*-N[\theta_1]-N[\theta_2]$, we note a unique 1-factor $F$, formed by edges $f_1=(-e_1-e_2+e_3,-e_2+e_3)$, $f_2=(e_1-e_2+e_3,2e_1-e_2+e_3)$, $f_3=(2e_1+e_3,2e_1+e_2+e_3)$,
$f_4=(2e_1+2e_2,2e_1+2e_2+e_3)$, $f_5=(2e_2+e_3,e_1+2e_2+e_3)$, $f_6=(2e_2-e_3,e_1+2e_2-e_3)$,
$f_7=(-e_1+2e_2-e_3,-e_1+e_2-e_3)$, etc. The copies of $K_2$ containing $f_1,\ldots,f_6$ can be taken dominated, by symmetry and forcedly, by the copies of $K_2$ containing $f_1-e_2$, $f_2+e_3$, $f_3+e_1$, $f_4+e_2$, $f_5+e_3$ and $f_6-e_3$ respectively. The 4-cycle induced in $S$ that contains $\theta_1$, also contains forcedly the vertices $\theta_1-e_1$, $\theta_1-e_1-e_3$ and $\theta_1-e_3$. But then, $f_7$ cannot be dominated in $S$, a contradiction.

Now, assume that the two corners are at distance 5, (Figure 6). They may be taken up to symmetry as $\theta_1=-e_1-e_2-e_3$ and $\theta_2=2e_1-e_2+e_3$. In $\Theta_0'-\Theta_0^*-N[\theta_1]-N[\theta_2]$ we observe a unique 1-factor $F$, formed by  edges $f_1=(-e_1+e_2-e_3,-e_1+2e_2-e_3)$, $f_2=(e_1-e_2-e_3,2e_1-e_2-e_3)$,
$f_3=(2e_2+e_3,e_1+2e_2+e_3)$, $f_4=(2e_1+2e_2,2e_1+2e_2-e_3)$, $f_5=(2e_1-e_3,2e_1+e_2-e_3)$, $f_6=(-e_1+e_3,-e_1+e_2+e_3)$, $f_7=(2e_2-e_3,e_1+2e_2-e_3)$, etc.
If the edge $(\theta_1,\theta_1-e_3)$ is in $S$, then
$f_1-e_1$, $f_6-e_3$, $f_3+e_2$, $f_4+e_1$ and $f_7-e_3$ dominate respectively $f_1$, $f_6$ $f_3$, $f_4$ and $f_7$. But then $f_5$ cannot be dominated in $S$, a contradiction. So, $F$
forces the 4-cycle with vertices $\theta_1$, $\theta_1-e_1$, $\theta_1-e_1-e_2$ and $\theta_1-e_2$ to be in $S$. In this case, the copies of $K_2$ associated to $f_1$, $f_2$, $f_7$ and $f_4$ are dominated respectively by the copies of $K_2$ containing $f_1-e_3$, $f_2-e_3$, $f_7+e_2$ and $f_4+e_1$. It follows that $f_5$ cannot be dominated by an edge at distance 1 from it in $\Lambda_3-\Theta_0'$, a contradiction.

It is easy to see that two corners at distance 6 or 8 do not allow even the definition of a 1-factor $F$ in $\Theta_0'-\Theta_0^*$ minus the two corners and their neighbors.

We pass to consider the different cases of four corners of $S$ in $\Theta_0'-\Theta_0^*$. The case of $S$ having three corners on the affine plane $<e_1,e_2>-e_3$ and one corner in the affine plane $<e_1,e_2>+e_3$, or viceversa, is readily seen to lead to no 1-factor $F$ in $\Theta_0'-\Theta_0^*$ minus these corners and their neighbors. Else, either:

{\bf Instance (A)}: If the four corners in $S$ are $\theta_1=-e_1-e_2-e_3$, $\theta_2=2e_1-e_2-e_3$, $\theta_3=-e_1+2e_2+e_3$ and $\theta_4=2e_1+2e_2+e_3$, then a 1-factor $F$ of $\Theta_0'-\Theta_0^*-\cup_{i=1}^4 N[\theta_i]$ is formed by the edges $f_1=(-e_1-e_2+e_3,-e_1+e_3)$,
$f_2=(-e_1+e_2-e_3,-e_1+2e_2-e_3)$,
$f_3=(-e_2+e_3,e_1-e_2+e_3)$,
$f_4=(2e_1-e_2+e_3,2e_1+e_3)$,
$f_5=(2e_1+e_2-e_3,2e_1+2e_2-e_3)$,
$f_6=(2e_2-e_3,e_1+2e_2-e_3)$.
We first rule out the case of
the edges $(\theta_1,\theta_1-e_1)$ and $(\theta_3,\theta_3+e_3)$ being in $S$ (or any other pair of edges in the same relative geometrical positions as these two, with respect to $\Theta_0'$). In this case, $f_1$ cannot be dominated by any copy of $K_2$: the two candidates, $f_1-e_1$ and $f_1+e_3$ cannot be in $S$. Because of this, three cases can be distinguished here up to symmetry, for the 4-cycles corresponding respectively to the four corners above, namely:

\begin{figure}[htp]\hspace*{18mm}%Figure 8
  \includegraphics[scale=0.45]{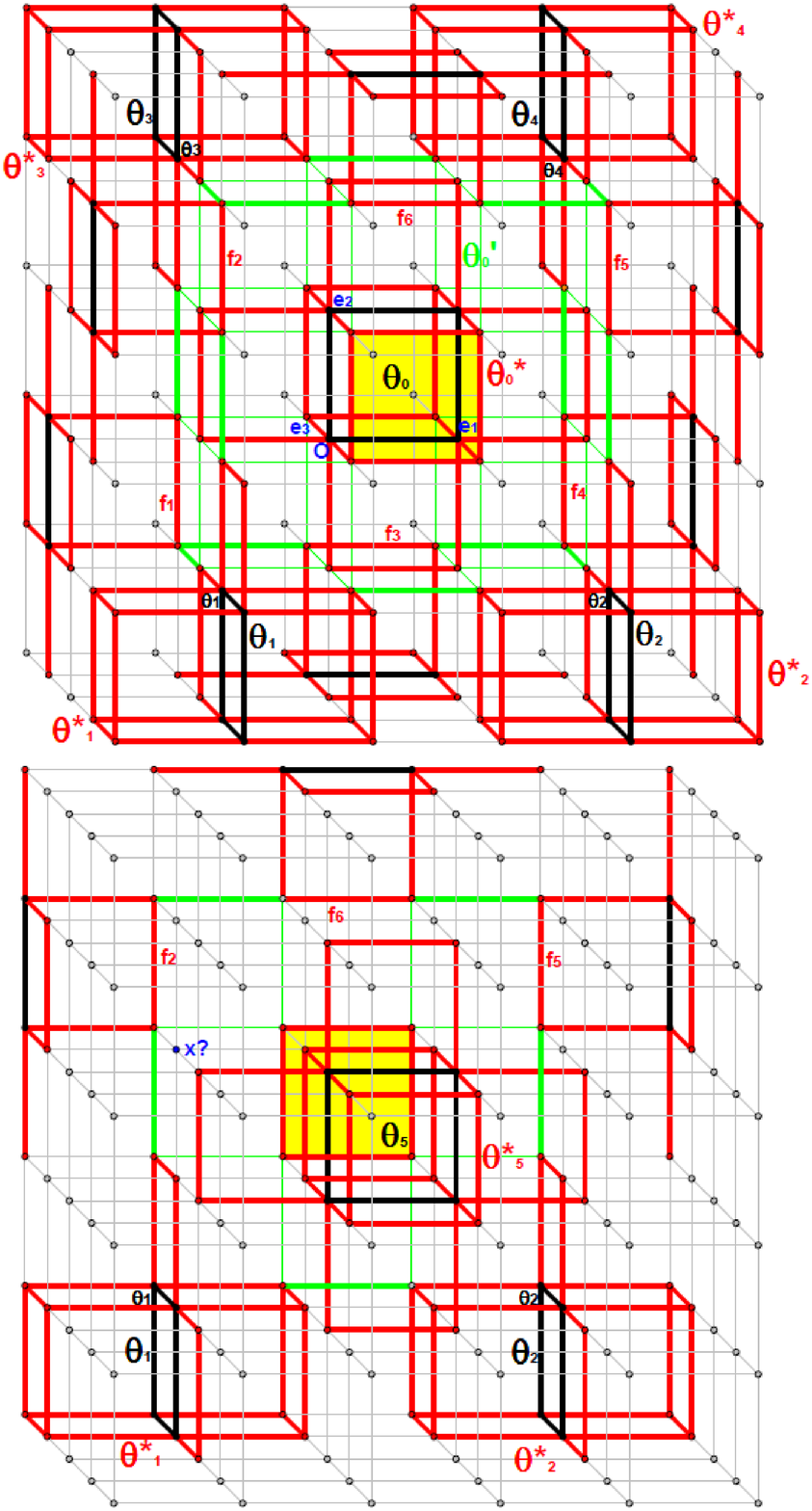}
\caption{Four corners: instance (A), case (b$_1$), subcase (b$_{11}$)}
\end{figure}

\begin{figure}[htp]\hspace*{18mm}%Figure 9
  \includegraphics[scale=0.45]{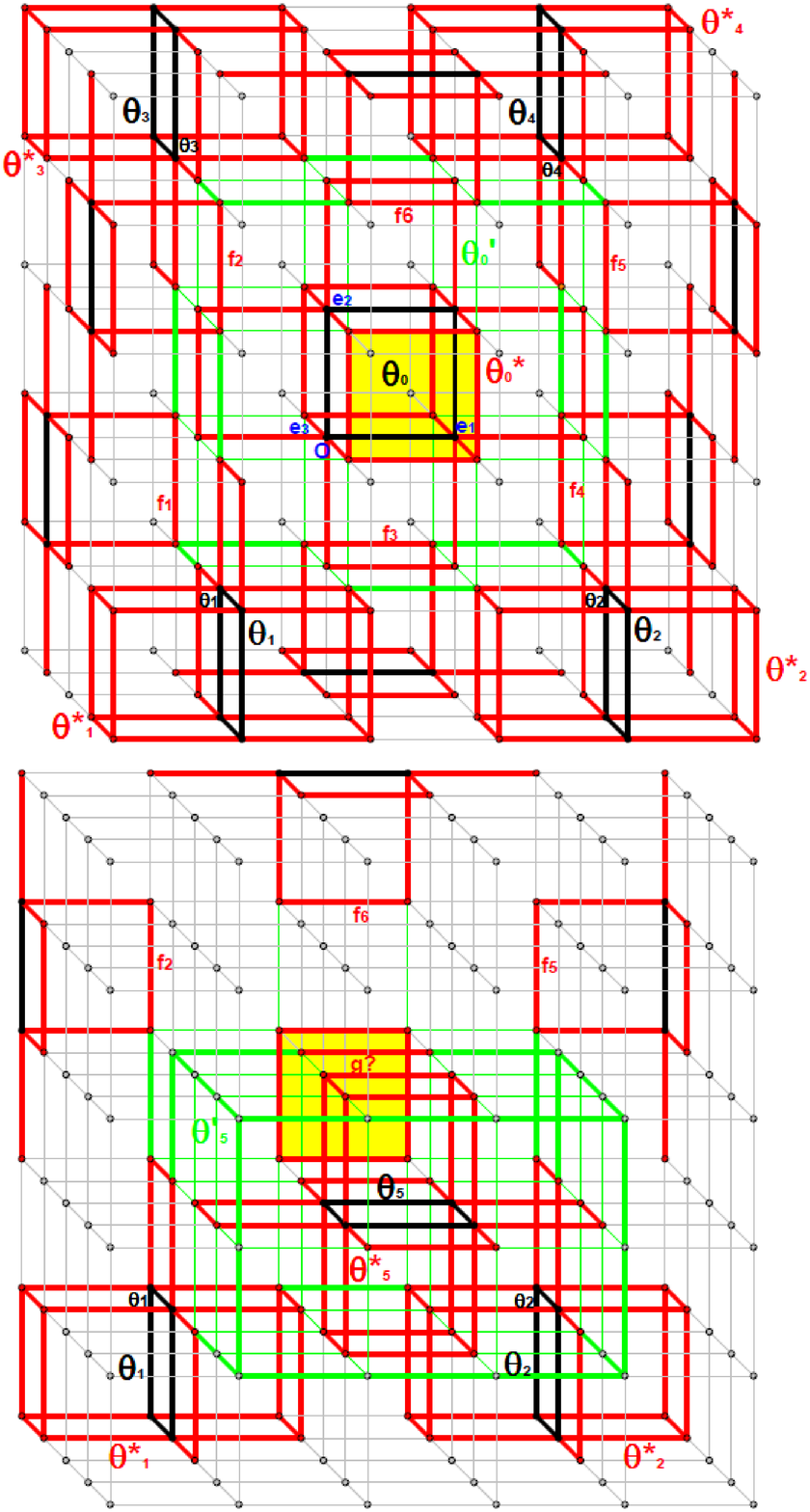}
\caption{Four corners, instance (A), case (b$_1$), subcase (b$_{12}$)}
\end{figure}

{\bf(a)} (Figure 7) $\Theta_1=(\theta_1,\theta_1-e_1,\theta_1-e_1-e_2,\theta_1-e_2)$,
$\Theta_2=(\theta_2,\theta_2+e_1,\theta_2+e_1-e_2,\theta_2-e_2)$,
$\Theta_3=(\theta_3,\theta_1-e_1,\theta_3-e_1-e_2,\theta_3-e_2)$,
$\Theta_4=(\theta_4,\theta_2+e_1,\theta_4+e_1-e_2,\theta_4-e_2)$.
Then the following edges must be in $S$, dominating forcedly the edges of $F$:
$f_1+e_3$, $f_2-e_3$, $f_3-e_2$, $f_4+e_3$, $f_5-e_3$, $f_6+e_2$.
The following 4-cycles are induced by $S$: $\Theta_5=(-e_1+e_2-e_3,-e_1+2e_2-e_3,-e_1+2e_2-2e_3,-e_1+e_2-2e_3)$ and
$\Theta_6=(2e_1+e_2-e_3,2e_1+2e_2-e_3,2e_1+2e_2-2e_3,2e_1+e_2-2e_3)$.
The graphs $\Theta'_5-\Theta^*_5$ and
$\Theta'_6-\Theta^*_6$ have the respective vertices $x=-3e_3$ and $y=e_1-3e_3$ as non-corner vertices, so they cannot dominate $z=-2e_3$ and $w=e_1-2e_3$, yielding a contradiction.

{\bf(b)} $\Theta_1=(\theta_1,\theta_1-e_2,\theta_1-e_2-e_3,\theta_1-e_3)$,
$\Theta_2=(\theta_2,\theta_2-e_2,\theta_2-e_2-e_3,\theta_2-e_3)$,
$\Theta_3=(\theta_3,\theta_3+e_2,\theta_3+e_2+e_3,\theta_3+e_3)$,
$\Theta_4=(\theta_4,\theta_4+e_2,\theta_4+e_2+e_3,\theta_4+e_3)$.
Then the following edges must be in $S$, dominating forcedly the edges of $F$:
$f_1-e_1$, $f_2-e_1$, $f_4+e_1$, $f_5+e_1$ and possibly:

{\bf(b$_1$)} (Figures 8-9) $f_6+e_2$, in which case: {\bf(b$_{11}$)} either $\Theta_5=\Theta_0-3e_3$ is in $S$ and dominates $\Theta_0-2e_3$, so that $x=-e_1+e_2-2e_3$ cannot be dominated by any of its neighbors; {\bf(b$_{12}$)} or $\Theta_5=(-3e_3,e_1-3e_3,e_1-4e_3,-4e_3)$ is in $S$, so the end-vertices of the edge $g=(e_2-2e_3,e_1+e_2-2e_3)$ cannot be dominated by $S$;

\begin{figure}[htp]\hspace*{18mm}%Figure 10
  \includegraphics[scale=0.45]{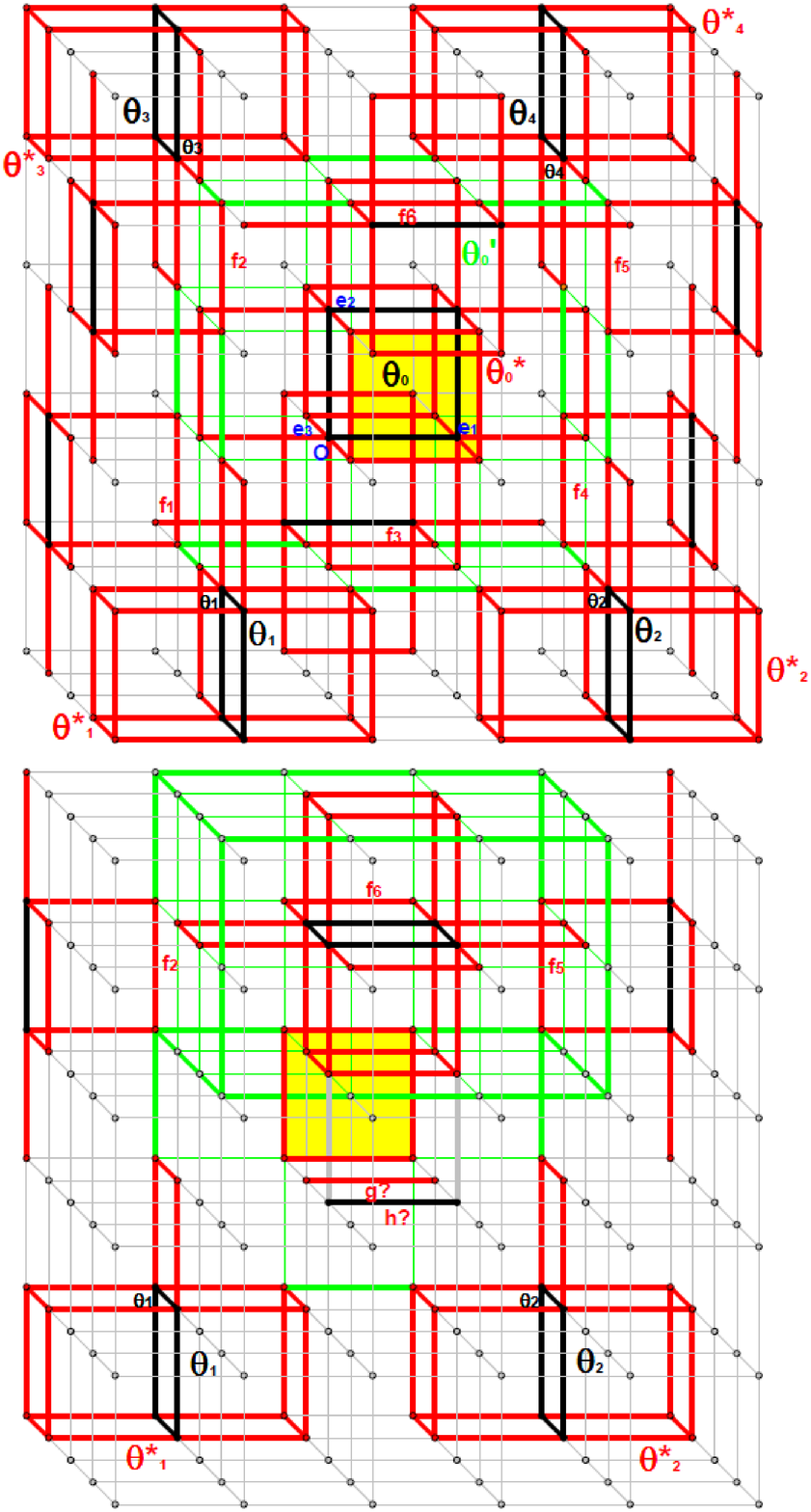}
\caption{Four corners, instance (A), case (b$_2$)}
\end{figure}

\begin{figure}[htp]\hspace*{18mm}%Figure 11
  \includegraphics[scale=0.45]{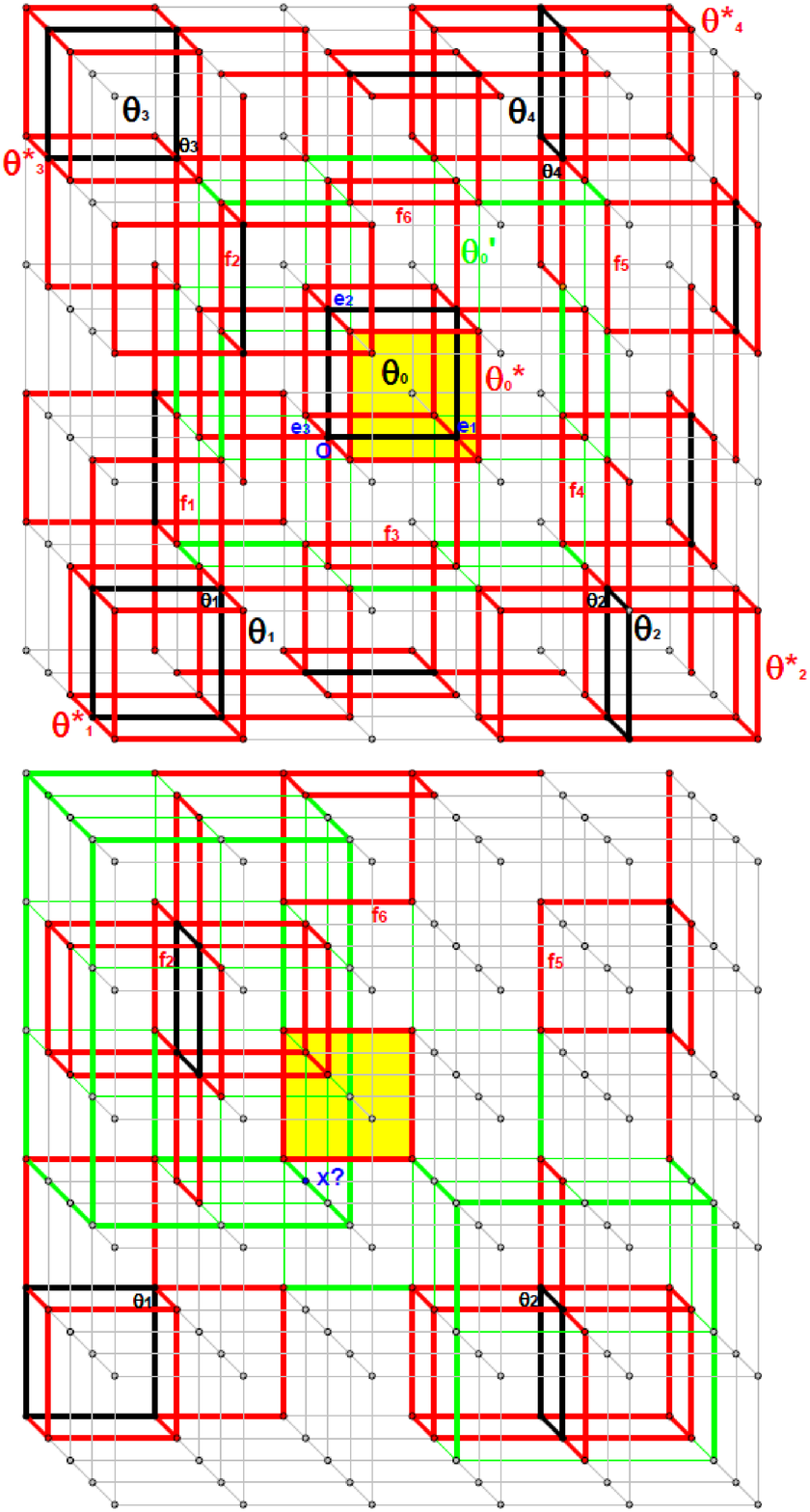}
\caption{Four corners, instance (A), case (c)}
\end{figure}

%\newpage

{\bf(b$_2$)} (Figure 10) $f_6-e_3$, in which case the end vertices of the edge $g=(e_2-2e_3,e_1+e_2-2e_3)$ cannot be in $S$ or dominated by $S$, since $h=(e_2-3e_3,e_1+e_2-3e_3)$ cannot be in $S$.

{\bf(c)} (Figure 11) $\Theta_1=(\theta_1,\theta_1-e_1,\theta_1-e_1-e_2,\theta_1-e_2)$,
$\Theta_2=(\theta_2,\theta_2-e_2,\theta_2-e_2-e_3,\theta_2-e_3)$,
$\Theta_3=(\theta_3,\theta_1-e_1,\theta_3-e_1-e_2,\theta_3-e_2)$,
$\Theta_4=(\theta_4,\theta_4+e_2,\theta_4+e_2+e_3,\theta_4+e_3)$.
Then the following edges must be in $S$, dominating forcedly the edges of $F$:
$f_1+e_3$, $f_2-e_3$, $f_3-e-2$, $f_4+e_1$, $f_5+e_1$, $f_6+e_2$. It follows that $x=-2e_3$ cannot be dominated by $S$.

Or {\bf Instance (B)}: For the rest, we need by symmetry only to consider the case in which the four corners of $S$ in $\Theta_0'-\Theta_0^*$ are $\theta_1=-e_1-e_2+e_3$, $\theta_2=2e_1-e_2+e_3$, $\theta_3=-e_1+2e_2+e_3$ and $\theta_4=2e_1+2e_2+e_3$. In the intersection of the affine plane $<e_1,e_2>-e_3$ and $\Theta_0'-\Theta_0^*$, a 1-factor $F$ is formed by the edges of the copies of $K_2$ that should be dominated externally (off $\Theta_0'$) by induced copies of $K_2$ in $S$ (parts themselves of 4-cycles induced by $S$). We may assume that this 1-factor is formed by the edges
$f_1=(-e_2-e_3,e_1-e_2-e_3)$, $f_2=(2e_2-e_3,e_1+2e_2-e_3)$,
$f_3=(-e_1-e_2-e_3,-e_1-e_3)$,
$f_4=(-e_1+e_2-e_3,-e_1+2e_2-e_3)$,
$f_5=(2e_1-e_2-e_3,2e_1-e_3)$ and
$f_6=(2e_1+e_2-e_3,2e_1+2e_2-e_3)$.
It is enough to consider by symmetry three cases of how $F$ could be dominated externally, as just mentioned, These cases have in common that
$f_1$ is dominated by $f_1-e_2$, $f_3$ by $f_3-e_3$, $f_4$ by $f_4-e_1$,
and differ in that:

{\bf(a)} (Figure 12)
$f_2$ is dominated by $f_2+e_2$, $f_5$ by $f_5+e_1$, $f_6$ by $f_6-e_3$;

{\bf(b)} (Figure 13)
$f_2$ is dominated by $f_2-e_3$, $f_5$ by $f_5-e_3$, $f_6$ by $f_6+e_1$;

{\bf(c)} (Figure 14)
$f_2$ is dominated by $f_2+e_2$, $f_5$ by $f_5-e_3$, $f_6$ by $f_6+e_1$.

In either case, by considering the dominating 4-cycle $\Theta_1=(-e_1-e_2-2e_3$, $-e_1-2e_3$, $-e_1-3e_3$, $-e_1-e_2-3e_3)$, the corresponding $\Theta_1'-\Theta_1^*$ contains two corners at distance 5, namely $x=-2e_2-e_3$ and $y=-2e_1+e_2-e_3$, which was ruled out above.

We just finished showing that there do not exist non-lattice like PDS$[Q_2]$\thinspace s in $\Lambda_3$. Thus, the only standing case of a PDS$[Q_2]$ in $\Lambda_3$ is the lattice-like one that remained by means of the commented programming code at the beginning of the present proof that leads to the generator matrix~(\ref{gm}) or its associated group epimorphism $\Phi:\mathbb{Z}^3\rightarrow \mathbb{Z}_{20}$. This establishes the statement of the theorem.
\end{proof}

\begin{figure}[htp]\hspace*{18mm}%Figure 12
  \includegraphics[scale=0.45]{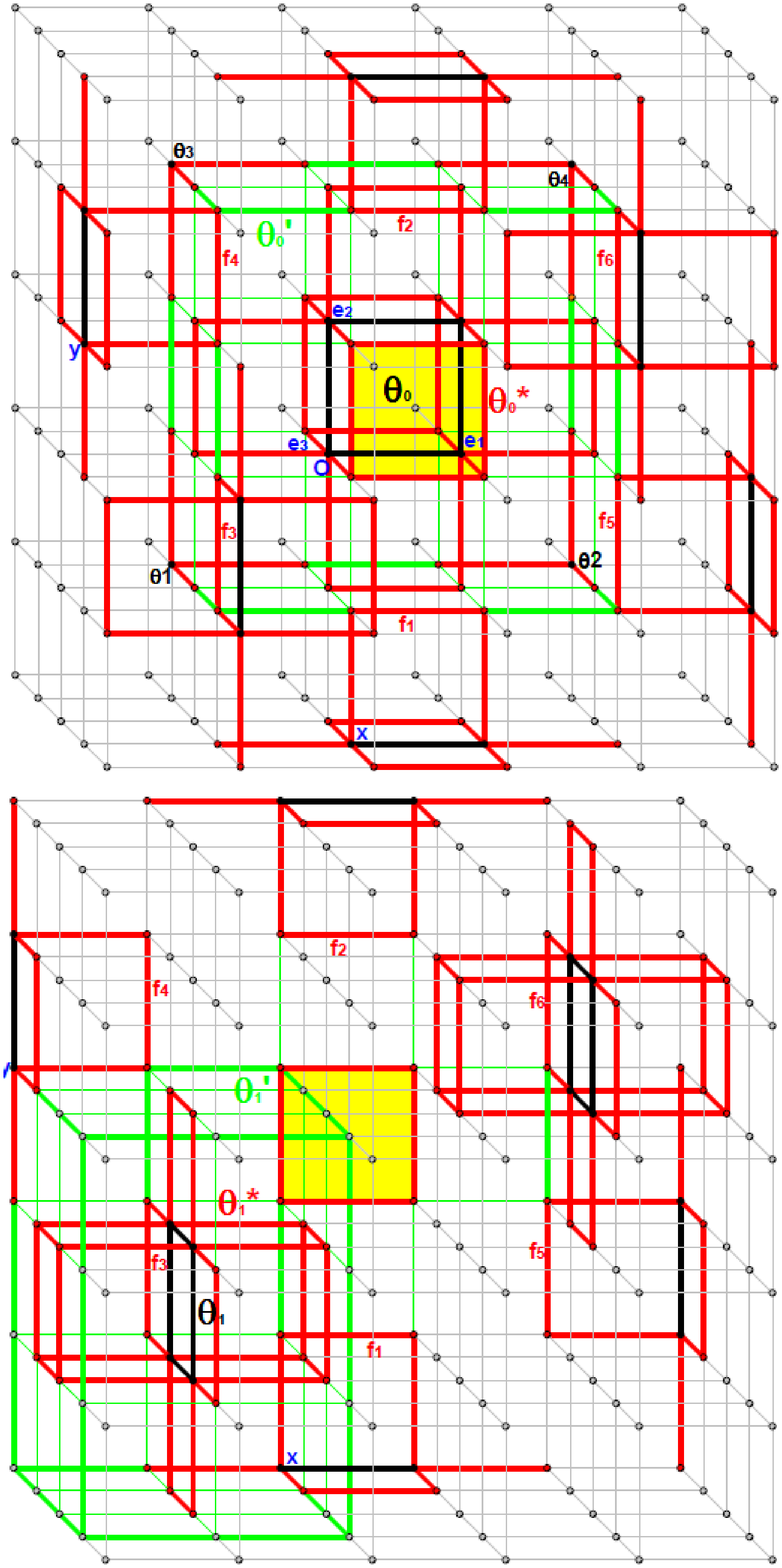}
\caption{Four corners, instance (B), case {\bf(a)}}
\end{figure}

\begin{figure}[htp]\hspace*{18mm}%Figure 13
  \includegraphics[scale=0.45]{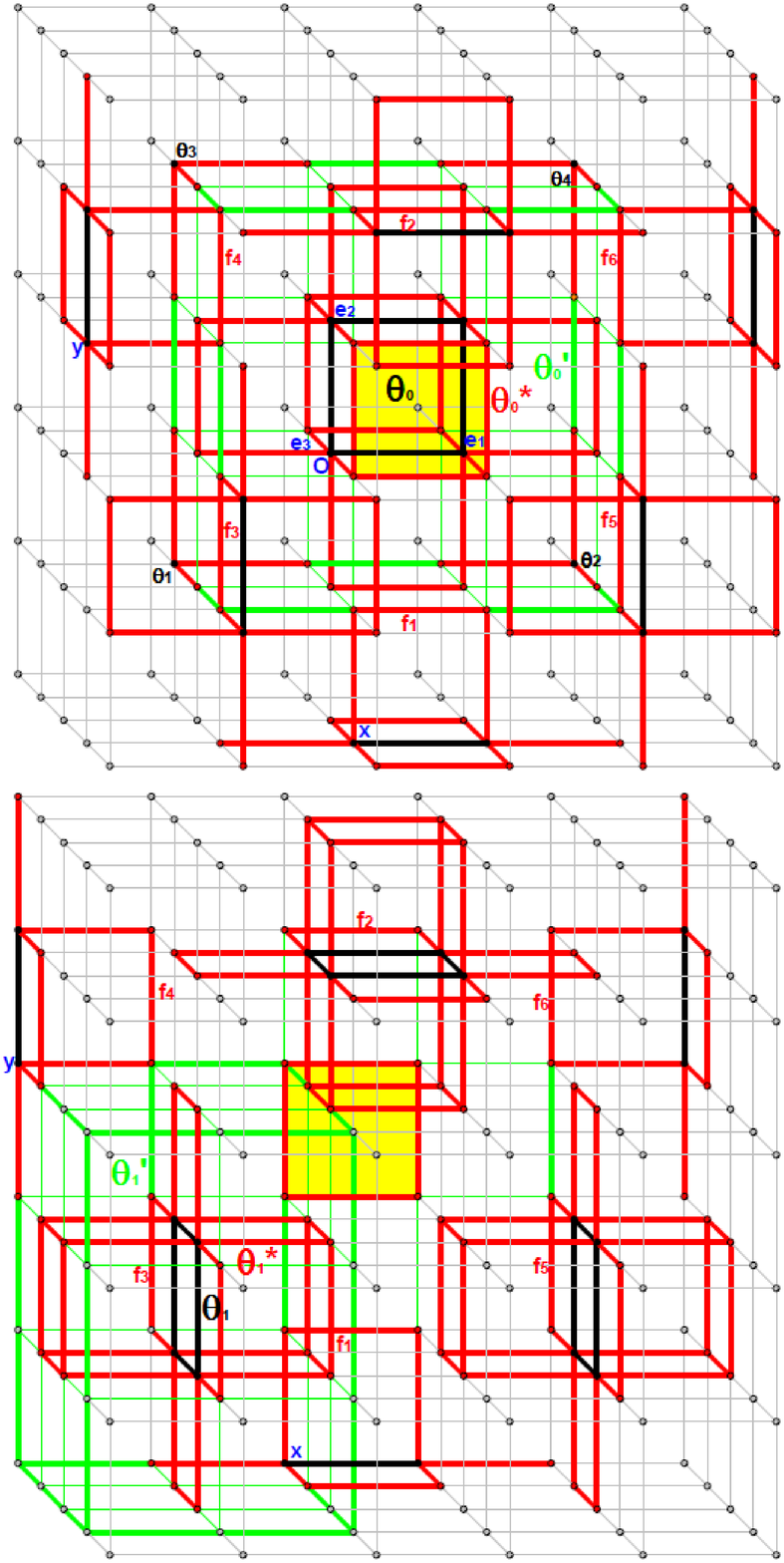}
\caption{Four corners, instance (B), case {\bf(b)}}
\end{figure}

\begin{figure}[htp]\hspace*{18mm}%Figure 14
  \includegraphics[scale=0.45]{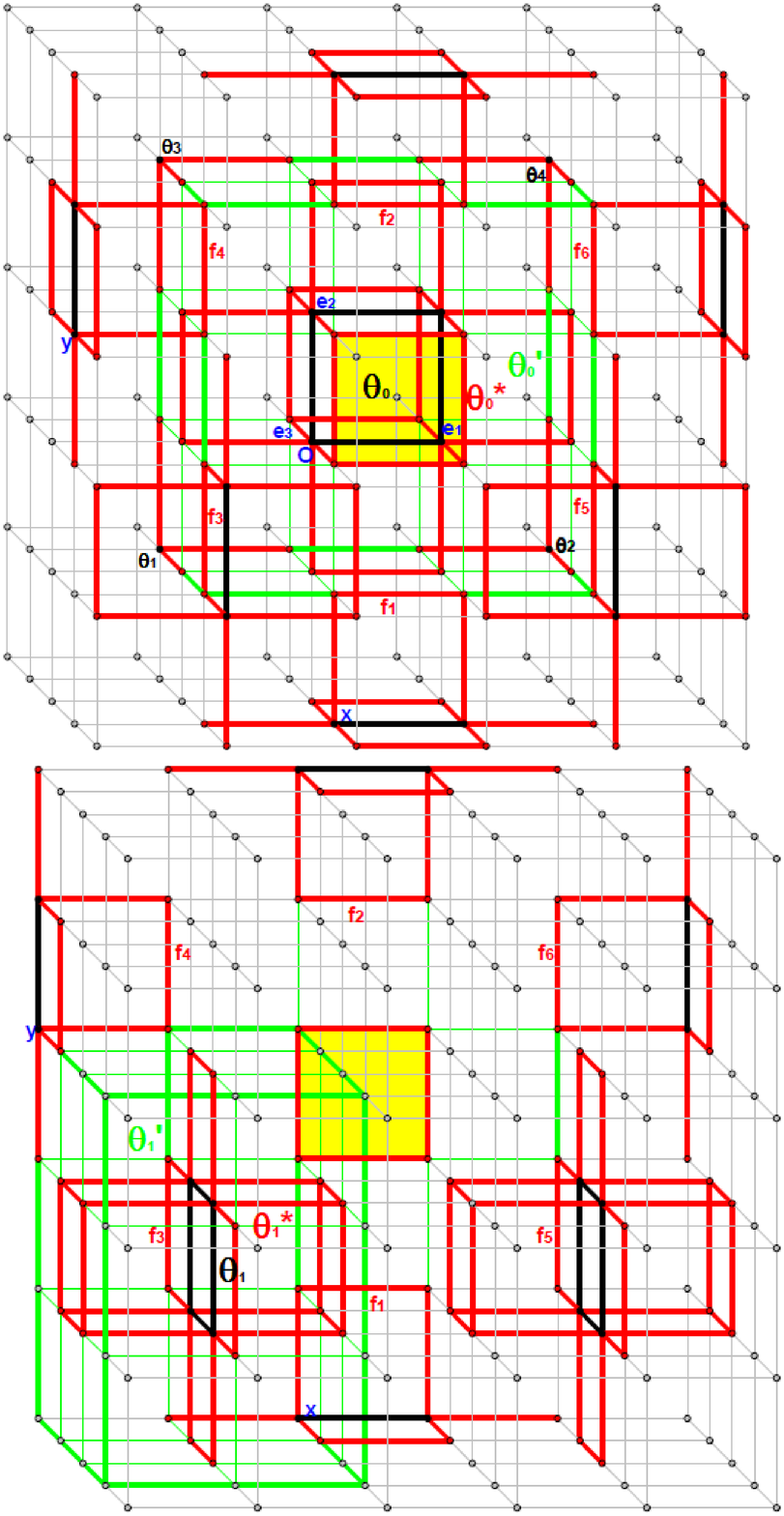}
\caption{Four corners, instance (B), case {\bf(c)}}
\end{figure}

\end{document}